# Data-Driven Planning of Plug-in Hybrid Electric Taxi Charging Stations in Urban Environments: A Case in the Central Area of Beijing

Huimiao Chen[1], Yinghao Jia[2], Zechun Hu[1], Guanglei Wu[1], Zuo-Jun Max Shen[2, 3, 4]

[1]Department of Electrical Engineering
Tsinghua University
Beijing, China
(chenhm15@mails.tsinghua.edu.cn)

[2]Department of Industrial Engineering
Tsinghua University
Beijing, China

[3]Department of Industrial Engineering and Operations Research
University of California, Berkeley
Berkeley, CA, USA

[4]Department of Civil and Environmental Engineering
University of California, Berkeley
Berkeley, CA, USA

*Abstract*—**Plug-in electric vehicles (PEVs) can contribute to the improvement of energy and environmental issues. Among different types of PEVs, plug-in hybrid electric taxis (PHETs) go in advance. In this study, we provide a spatial and temporal PHET charging demand forecasting method based on one-month global positioning system (GPS)-based taxi travel data in Beijing. Then, using the charging demand forecasting results, a mixed integer linear programming (MILP) model is formulated to plan PHET charging stations in the central area of Beijing. The model minimizes both investment and operation costs of all the PHET charging stations and takes into account the service radius of charging stations, charging demand satisfaction and rational occupation rates of chargers. At last, the test of the planning method is carried out numerically through simulations and the analysis is complemented according to the results.**

## I. Introduction

As a cleaner mode of transport, plug-in electric vehicles (PEVs) have been long regarded as a promising tool to reduce traffic emissions and petroleum dependence. In recent years, grand development plans with the goal of proliferating PEVs have been brought forward around the world [1]. For example, in China, the government is aiming at an overall ownership of 5 million new energy vehicles by 2020 [2].

Although great efforts have been made to electrify transportation, barriers to large-scale commercialization of PEVs exist in various aspects. Among others, a major hurdle is the so-called range anxiety phenomenon, which is associated with the PEV user's fear of running out of battery power en-route and keeps customers away from PEVs. Therefore, it is vitally important to solve the lack of enabling charging infrastructure.

Driven by the urgent needs of satisfying the increasing PEV charging demands, the planning of PEV charging facilities has been a research hotspot for years. A widely adopted approach to the planning problem seeks to locate charging facilities near the activity centers of PEV users [3]. But the planning objective and method vary. Reference [4] firstly uses a two-step screening method with environmental factors and service radius of PEV charging stations considered to identify the candidate PEV charging stations, and then optimizes the total costs by a modified primal-dual interior point algorithm. In [5], for the maximization of the social welfare associated with the coupled power and transportation networks, the allocation problem of a given number of public charging stations among metropolitan areas is formulated as a mathematical program with complementarity constraints, and is solved by an active-set algorithm. Authors of [6] develop a maximal covering model to locate a fixed number of charging stations in central urban areas. In [7], authors develop a multi-objective PEV charging station planning method which can ensure charging service while reducing energy losses and voltage fluctuations of distribution networks. Reference [8] is devoted to addressing the collaborative planning strategy for integrated power distribution and EV charging systems.

However, the above research is done based on assumptive behaviors of PEV users. Regarding literature considering the real-world travel data of vehicles, Cavadas *et al.* [9] take into account the possibility of transferring charging demands in the successive stops of the same trip when planning PEV charging stations with the purpose to maximize the satisfied charging demands under a fixed budget, and test their approach based on vehicle trip data from a survey questioned a sample of 10,000 participants; reference [10] studies PEV charger location problems, analyzes the impact of public charging infrastructure deployment on increasing electric miles traveled, and carries out the case studies using the global positioning system (GPS) based travel survey data collected in the greater Seattle metropolitan area; according to the vehicle activity data obtained from a GPS tracked household travel survey in Austin, Texas, reference [11] examines the role of public charging infrastructure in increasing the share of driving on

This work was supported in part by the National Key Research and Development Program (2016YFB0900103) and the National Natural Science Foundation of China (51477082).

electricity of plug-in hybrid electric vehicles (PHEVs); Cai *et al.* [12] evaluate how travel patterns mined from big data can inform public charging infrastructure development, and assess environmental impacts for different charging infrastructure siting scenarios.

Different from the above papers, in this text, we focus on the planning of plug-in hybrid electric vehicle (PHET) charging stations in urban environments and aim to obtain the optimal scheme mathematically based on a large number of real-world taxi travel data. This work is driven by 1) taxi is a public means of transportation and is expected to contribute significantly to transportation electrification; 2) PHET play an important role at the initial development stage of electric vehicles and thereby necessitates the corresponding charging station planning. The main procedures and contributions of the paper are summarized below.

- Process the raw data, and select the taxi dwelling with an appropriate duration threshold to separate the whole travels and to estimate the temporal distribution of PHET charging demands.

- Divide the central area of Beijing into a number of equal square sectors and utilize the data to forecast the spatial distribution of PHET charging demands based on that.

- Formulate a mixed integer linear programming (MILP) model to plan PHET charging stations, which seeks for the minimum investment and operation costs of all the PHET charging stations and takes into account the service radius of charging stations, charging demand satisfaction and rational occupation rate of chargers.

The remainder of the paper is organized as follows. Section II describes the GPS-based taxi travel data in Beijing and the processing method. In Section III, the data-driven spatial and temporal PHET charging demand forecasting method is presented. In Section IV, the mathematical model of PHET charging station planning is formulated. Then, Section V shows the simulation results and the analysis. Finally, conclusions are drawn in Section VI.

## II. DESCRIPTION AND PROCESSING OF GPS-BASED TAXI TRAVEL DATA

The data set used in this paper is GPS-based taxi travel records of 29709 taxis in Beijing (approximately $\pi = 44.3\%$ of the fleet[1]) from July 1st to July 31st in 2016, which were collected by cellphone or on-board device.

In order to clean up the raw data, we remove the duplicated and incorrect data points. And after that, the data set includes a total of 2.65 billion data points, which covers over 7 million trips. Each data point contains the timestamp up to seconds, the taxi ID, and the taxi location (in longitude and latitude). In Table I, a sample of the taxi travel records in the data set is presented.

[1] Currently, the total number of taxis in Beijing is approximately 67000 [13].

TABLE I. A SAMPLE OF TAXI TRAVEL RECORDS

| Taxi ID | Timestamp | Longitude | Latitude |
|---|---|---|---|
| 1140 | 20160701053632 | 116.467003 | 39.928562 |

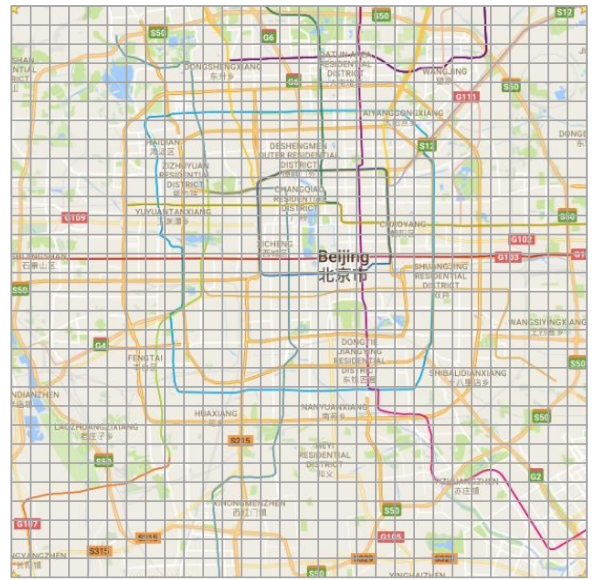

Fig. 1. The schematic illustration of dividing the central area of Beijing into 30×30 square sectors.

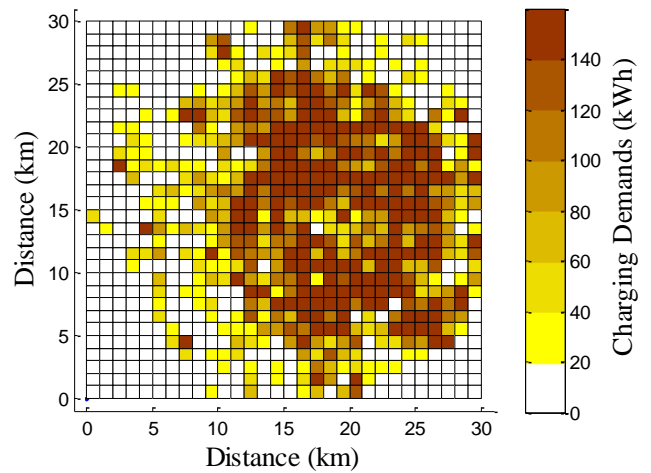

Fig. 2. The spatial distribution of PHET charging demands from 2 pm to 3 pm, July 11st, 2016.

In this research, we focus on the charging station planning for PHETs, i.e., the taxis which are electricity or petroleum powered. And we suppose that the traditional taxis are replaced by PHETs, and the PHETs give priority to consuming the battery power. Due to hybrid power, we assume that no range anxiety phenomenon occurs, i.e., the travel behaviors of the taxis remain unchanged after the adoption of PHETs [10], [12], [14]. Additionally, we regard the dwelling time of at least 30 minutes as the available recharging time windows[2], when charging demands may arise. Also, dwelling with no less than 30 minutes is used to separate the whole travels of a taxi into a series of segments without recharging. Then, the charging demands are temporally distributed during the valid dwelling. Further, we can extract all the location information of valid dwelling and the corresponding vehicle miles traveled (VMT) from the previous valid dwelling to this valid dwelling (VMT-D), from the data set. In Section III, the processed data are used for PHET charging demand forecasting.

III. DATA-DRIVEN SPATIAL AND TEMPORAL PHET CHARGING DEMAND FORECASTING

For the sake of forecasting the spatial and temporal PHET charging demands, we discretely divide a day into 24 time slots (the duration of a time slot is 1h, i.e., $\Delta t = 1h$) the central area of Beijing into 30×30 sectors (each sector is a 1km×1km square, see Fig. 1). Let $\mathbb{U}$ and $\mathbb{T}$ denote the sets of the sectors and the time slots in the concerned period, respectively.

Based on the processed data in Section II, given a time slot $t$, we can obtain the distribution of the number of PHETs with charging demands (valid dwelling) in all the sectors[3], which can be expressed by a number distribution matrix $\mathbf{N}^t$ with dimensions 30×30. Let $l_{x,y,i}^t, i=1,\cdots,\mathrm{N}_{x,y}^t$ ($\mathrm{N}_{x,y}^t$ is the element of $\mathbf{N}^t$) denotes the VMT-D of the PHET $i$ in sector $(x,y) \in \mathbb{U}$ at time slot $t$. Then, for any time slot $t$, we can obtain a charging demand distribution matrix $\mathbf{D}^t$, of which elements $\mathrm{D}_{x,y}^t$ can be calculated by (1).

$$\mathrm{D}_{x,y}^t = \sum_{i=1}^{\mathrm{N}_{x,y}^t} \xi \cdot \min\left(l_{x,y,i}^t, d\right), \forall (x,y) \in \mathbb{U} \quad (1)$$

where $\xi$ (=10kWh/50km) is the battery power consumption per kilometer in all-electric mode and $d$ (=50km) is the all-electric range of PHET. Note that for a PHET, its charging demands calculated by the VMT-D probably are more or less than its actual charging demands, but for a number of PHETs, the variations of individual PHET charging demands tend to cancel out each other. Fig. 2 shows the PHET charging demand distribution from 2 pm to 3 pm, July 11st, 2016.

All the $\mathbf{D}^t$, $t \in \mathbb{T}$ include the spatial and temporal distribution information of PHET charging demands.

IV. MATHEMATICAL MODEL FOR PLANNING OF PHET CHARGING STATIONS

In this section, the mathematical model is formulated to plan the PHET charging stations based on the data-driven spatial and temporal charging demands forecasting. The model aims to achieve the minimum costs for PHET charging stations, including both the investment costs and the operation costs. The objective function of the model is expressed as below.

$$\min \sum_{(i,j)\in\mathbb{U}} \left[ I_{i,j}\left(K_{i,j} \cdot \left(C^{\mathrm{ch}}+C^{\mathrm{ca}}P\right)+C_{i,j}^{\mathrm{la}}A_{i,j}^{\mathrm{st}}+C^{\mathrm{ot}}\right) + I_{i,j}C_{i,j}^{\mathrm{op}} \right] \quad (2)$$

$$= \sum_{(i,j)\in\mathbb{U}} \begin{bmatrix} I_{i,j}\left(K_{i,j}\cdot\left(C^{\mathrm{ch}}+C^{\mathrm{ca}}P\right)\right) \\ +C_{i,j}^{\mathrm{la}}\left(A^{\mathrm{st,fi}}+A^{\mathrm{st,va}}S_{i,j}\right)+C^{\mathrm{ot}}\right) \\ +I_{i,j}\left(C^{\mathrm{op,fi}}+C^{\mathrm{op,va,1}}\sum_{t\in\mathbb{T}}E_{i,j}^t+C^{\mathrm{op,va,2}}S_{i,j}\right) \end{bmatrix} \quad (3)$$

$$= \sum_{(i,j)\in\mathbb{U}} \begin{bmatrix} K_{i,j}\cdot\left(C^{\mathrm{ch}}+C_{i,j}^{\mathrm{la}}A^{\mathrm{st,va}}+C^{\mathrm{op,va,2}}+C^{\mathrm{ca}}P\right) \\ +I_{i,j}\left(C^{\mathrm{op,fi}}+C_{i,j}^{\mathrm{la}}A^{\mathrm{st,fi}}+C^{\mathrm{ot}}+C^{\mathrm{op,va,1}}\sum_{t\in\mathbb{T}}E_{i,j}^t\right) \end{bmatrix} \quad (4)$$

where $I_{i,j}$ and $K_{i,j}$ are decision variables, and $I_{i,j}$ denotes whether to set up a charging station in sector $(i,j)$ (binary variable, 1: building a charging station; 0: no charging station) and $K_{i,j}$ denotes the number of chargers (integer variable); $C^{\mathrm{ch}}$ is the cost of a charger; $C^{\mathrm{ca}}$ is the capacity charge, i.e., the cost per kWh of the installed capacity connected to the power grids; $P$ (=10kW) is the rated charging power of the charger; $C_{i,j}^{\mathrm{la}}$ is the land costs per unit area in sector $(i,j)$; $A_{i,j}^{\mathrm{st}}$ is the area of the charging station in sector $(i,j)$; $C^{\mathrm{ot}}$ is the other investment costs; $C_{i,j}^{\mathrm{op}}$ is the operation costs; $A^{\mathrm{st,fi}}$ and $A^{\mathrm{st,va}}$ are the fixed area of a charging station and the coefficient of the variable area of a charging station with respect to the number of chargers $K_{i,j}$, respectively; $C^{\mathrm{op,fi}}$ is the fixed operation costs of a charging station; $C^{\mathrm{op,va,1}}$ is the coefficient of the first part of variable operation costs, which are proportional to the total output energy $\sum_{t\in\mathbb{T}}E_{i,j}^t$; $C^{\mathrm{op,va,2}}$ is the coefficient of the second part of variable operation costs, which are proportional to the number of chargers $K_{i,j}$. In (2) and (3), the first and second terms are the investment costs and the operation costs, respectively. And in (4), the expression is converted into the form of $\alpha \cdot K_{i,j} + \beta \cdot I_{i,j}$ ($\alpha$ and $\beta$ are coefficients).

In consideration of that the total output energy of charging stations is nearly equal to the charging demands of the PHET fleet, expressed in (5), the objective function (4) can be further rewritten in a linear form as (6) by removing the terms including $E_{i,j}^t$.

$$\sum_{(i,j)\in\mathbb{U}}\sum_{t\in\mathbb{T}}E_{i,j}^t \approx \sum_{(i,j)\in\mathbb{U}}\sum_{t\in\mathbb{T}}\mathrm{D}_{x,y}^t \approx \text{constant} \quad (5)$$

$$\min \sum_{(i,j)\in\mathbb{U}} \begin{bmatrix} K_{i,j}\cdot\left(C^{\mathrm{ch}}+C_{i,j}^{\mathrm{la}}A^{\mathrm{st,va}}+C^{\mathrm{op,va,2}}+C^{\mathrm{ca}}P\right) \\ +I_{i,j}\left(C^{\mathrm{op,fi}}+C_{i,j}^{\mathrm{la}}A^{\mathrm{st,fi}}+C^{\mathrm{ot}}\right) \end{bmatrix} \quad (6)$$

---

[2] According to the parameters of some current plug-in hybrid electric vehicles on the market [15], [16], herein, we suppose the rated charging power is 10 kW, and let PHETs' battery capacity be 10 kWh and their all-electric range be 50 km. Derived from these parameters, the dwelling time no less than 30 minutes are approximately set as valid periods of recharging.
[3] We count the number of taxis whose beginning time of the valid dwelling time is within the time slot.

The constraints are presented as follows.

$$0 \leq I_{i,j} \leq \bar{I}_{i,j}, \forall (i,j) \in \mathbb{U} \quad (7)$$

$$I_{i,j} \leq K_{i,j} \leq \bar{K}_{i,j} \cdot I_{i,j}, \forall (i,j) \in \mathbb{U} \quad (8)$$

$$0 \leq S_{i,j,m,n}^t \leq \bar{S}_{i,j,m,n}^t, \forall (i,j),(m,n) \in \mathbb{U}, t \in \mathbb{T} \quad (9)$$

$$\sum_{(i,j) \in \mathbb{U}} S_{i,j,m,n}^t \cdot P \cdot \Delta t = \gamma \cdot \tilde{D}_{m,n}^t,$$

$$\forall (m,n) \in \mathbb{U}, t \in \mathbb{T} \quad (10)$$

$$\sum_{(m,n) \in \mathbb{U}} S_{i,j,m,n}^t \leq \eta \cdot K_{i,j}, \forall (i,j) \in \mathbb{U}, t \in \mathbb{T} \quad (11)$$

where $\bar{I}_{i,j}$ is binary input parameter of candidate sectors of PHET charging stations (1: sector $(i,j)$ is a candidate sector of PHET charging stations; 0: otherwise); $\bar{K}_{i,j}$ is the upper bound of the number of chargers which can be installed in sector $(i,j)$; $S_{i,j,m,n}^t$ denotes the number of chargers in sector $(i,j)$, which are used to satisfy the charging demands from sector $(m,n)$, at time slot $t$ (i.e., $S_{i,j,m,n}^t$ chargers in sector $(i,j)$ serve the PHETs coming from sector $(m,n)$ at time slot $t$, and (10) provides the mathematical explanation to $S_{i,j,m,n}^t$); $\bar{S}_{i,j,m,n}^t$ is the upper bound of $S_{i,j,m,n}^t$, which can be calculated by (12); $\gamma$ is the PHET penetration rate, i.e., the proportion of PHETs to the taxi fleet; $\tilde{D}_{m,n}^t$ is the charging demands of PHETs in sector $(i,j)$ at time slot $t$ if all the taxis are replaced by PHETs, calculated by (13); $\eta$ is the upper bound of the occupancy rate of chargers.

$$\begin{cases} \bar{S}_{i,j,m,n}^t = \Gamma, (|i-m|+|j-n|) \cdot a \leq r \\ \bar{S}_{i,j,m,n}^t = 0, (|i-m|+|j-n|) \cdot a > r \end{cases},$$

$$\forall (i,j),(m,n) \in \mathbb{U}, t \in \mathbb{T} \quad (12)$$

$$\tilde{D}_{m,n}^t = D_{m,n}^t / \pi, \forall (m,n) \in \mathbb{U}, t \in \mathbb{T} \quad (13)$$

where $\Gamma$ is a sufficiently large positive value; $a$ ($=1\text{km}$) is the length of the side of a square sector; $r$ is the service radius of a charging station.

In the above constraints, (7)-(9) are constraints of variables ranges; (10) are balance constraints of the output energy and the charging demands; (11) describe the limitation of the outputs of charging stations. Note that: 1) $S_{i,j,m,n}^t$ is not an integer variable because it is actually the average number of chargers within time slot $t$; 2) a PHET is supposed to go to a charging stations for recharging if their distance is within the predefined service radius (see (9) and (12)); 3) we calculate the Manhattan distance between a PHET and a charging station (see (12)) because this paper focuses the case in Beijing where the roads are typically vertical and horizontal; 4) a higher occupancy rate of chargers reduces the number of chargers, and accordingly brings lower investment, but longer average working time of chargers and longer waiting time of PHETs.

Thus, setting an upper bound of occupancy rates is necessary (see (11)); 5) the initial investment costs in the model, e.g.,

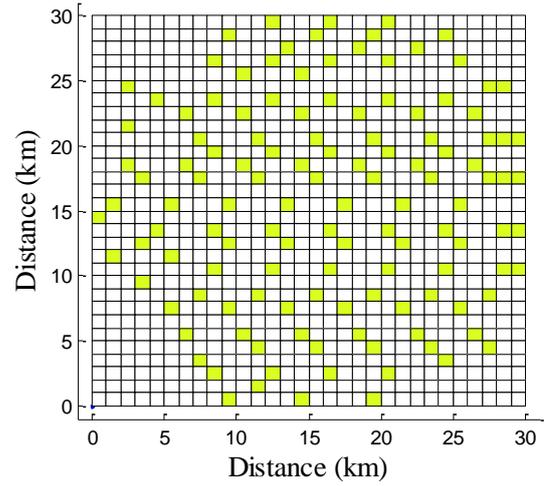

Fig. 3. The distribution of candidate sectors of PHET charging stations (green sectors are candidate sectors).

TABLE II. PARAMETER SETTINGS

| | |
|---|---|
| $C^{ch}$ | 1140 RMB/(month each charger) (equivalent monthly value) |
| $C_{i,j}^{la}$ | [300 RMB/month, 400 RMB/month] (linearly decreasing from the center to the edge in the map shown in Fig. 1) |
| $C^{op,va,2}$ | 50 RMB/(month each charger) |
| $C^{op,fi}$ | 500 RMB/(month each station) |
| $C^{ca}$ | 70 RMB/(kW month) |
| $C^{ot}$ | 912 RMB/(month each station) |
| $A^{st,va}$ | 20 m²/each charger |
| $A^{st,fi}$ | 100 m²/each station |
| $\bar{K}_{i,j}$ | all set as 150 chargers |
| $P$ | 10kW/each charger |
| $\gamma$ | 40% and 50% |
| $\eta$ | 80% |
| $r$ | 4km and 6km |

$K_{i,j} \cdot C^{ch}$, need to be converted into the equivalent costs in the concerned period (one month in this paper) leveraging the capital recovery factor [17].

The formulated model (6)-(11) of PHET charging station planning is a MILP problem, which can be solved efficiently via many available commercial optimizers.

## V. CASE STUDIES

### A. Parameter Settings

As aforementioned, we study the PHET charging station planning in the central area of Beijing utilizing the GPS-based taxi travel data from July 1st to July 31st in 2016. Same as that in charging demand forecasting, the map is divided into 30×30 square sectors and time is discretized hourly. The candidate sector distribution ($\bar{I}_{i,j}$) is shown in Fig. 3. All the candidate sectors are selected from the sectors with positive charging demands with the condition that the union service areas of the candidates cover all the sectors with charging demands. The other parameters used in the case are listed in Table II.

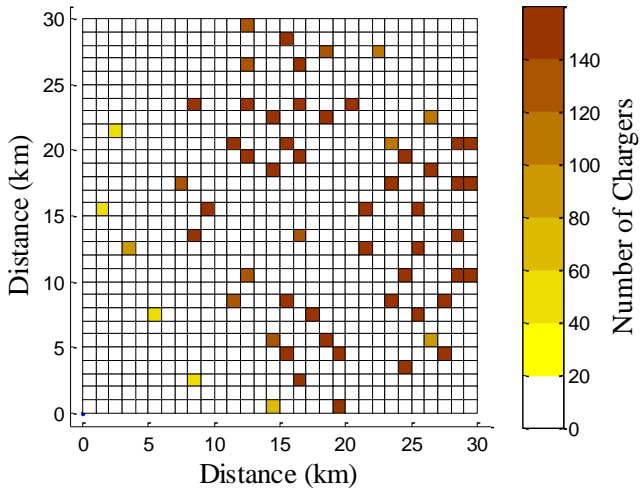

Fig. 4. The charger distribution results (PHET penetration rate 40%, service radius 4km).

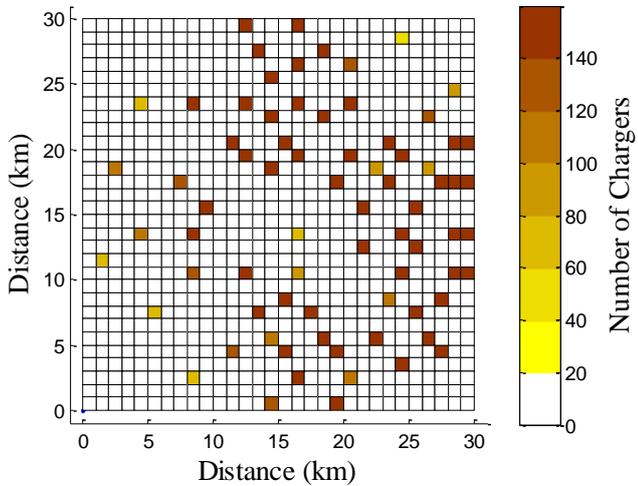

Fig. 5. The charger distribution results (PHET penetration rate 50%, service radius 4km).

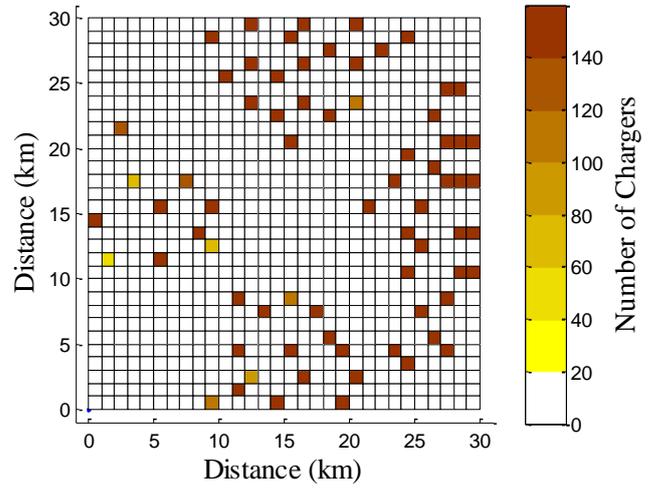

Fig. 6. The charger distribution results (PHET penetration rate 50%, service radius 6km).

*B. Results and Analysis*

We carry out the simulations when the PHET penetration rate and the service radius are (40%, 4km) (50%, 4km) and (50%, 6km), respectively. The charger distribution results are shown in Figs. 4, 5 and 6, respectively. And the corresponding equivalent monthly costs, including the investment and operation costs, are respectively 66,358,072 RMB and 82,971,340 RMB and 80,877,292 RMB. The simulations run for about 2 hours each time.

When the service radius is 4km, the cost results show that the payment is approximately proportional to the number of PHETs. And when the PHET penetration rate is 50%, the larger service radius brings lower costs. Also, according to the charger distribution results (see Figs 5 and 6), we can observe that the charging station number decreases as the service radius becomes larger. However, a larger service radius actually intensifies range anxiety phenomenon and transfers the partial economic loss to the PHET drivers. Thus, a rational service radius should be determined before the planning.

Comparing the candidate sector distribution and the charger distributions, it can be seen that the charging stations tend to be away from the center of the city because the land fee declines from the center to the edge. Additionally, the chargers tend to be centralized due to the fixed costs of charging stations.

The planning model was programmed in Matlab 2014a and solved by an embedded Cplex solver [16] with the YALMIP interface [17]. We implemented all the above numerical tests in a computational environment with Intel(R) Core(TM) i5-3210M CPUs running at 2.50 GHz with 8 GB random access memory.

## VI. Conslusions

In this paper, we present a planning method for PHETs in urban environments based on the GPS-based taxi travel data from Beijing. There are two main steps: 1) utilizing the data to forecast the spatial and temporal PHET charging demand distribution; 2) formulating an MILP model to plan PHET charging stations, which seeks for the minimum investment and operation costs of all the PHET charging stations with the constraints of the service radius of charging stations, charging demand satisfaction and rational occupation rate of chargers.

Case studies provide the analysis of the planning results and the conclusions include: 1) charging stations tend to be away from the center of the city; 2) the charging station number should be as fewer as possible; 3) too large service radius is not suggested even though it can reduce the total costs. In future work, we will consider the heterogeneous PHET population and further plan the charging facilities for pure electric taxis.